\documentclass[12pt]{article}

\usepackage{natbib}
\usepackage{amsthm,amsfonts,amssymb}
\usepackage{amsmath,thmtools}
\usepackage{setspace,graphicx,enumitem,soul,url}
\usepackage{tikz-cd}
\usepackage{authblk}

\newtheorem*{example}{Example}

\usepackage{authblk}
\author[1]{Martina Patone}
\author[1,2,3]{Li-Chun Zhang}
\affil[1]{\em \small University of Southampton, UK (email: L.Zhang@soton.ac.uk)}
\affil[2]{\em \small Statistics Norway, Norway}
\affil[2]{\em \small University of Oslo, Norway}
\title{Incidence weighting estimation\\ under bipartite incidence graph sampling}
\date{}

\begin{document}

\maketitle

\begin{abstract}
\noindent
Bipartite incidence graph sampling provides a unified representation of many sampling situations for the purpose of estimation, including the existing unconventional sampling methods, such as indirect, network or adaptive cluster sampling, which are not originally described as graph problems. We develop a large class of linear estimators based on the edges in the sample bipartite incidence graph, subjected to a general condition of design unbiasedness. The class contains as special cases the classic Horvitz-Thompson estimator, as well as the other unbiased estimators in the literature of unconventional sampling, which can be traced back to Birnbaum and Sirken (1965). Our generalisation allows one to devise other unbiased estimators, thereby providing a potential of efficiency gains in applications. Illustrations are given for adaptive cluster sampling, line-intercept sampling and simulated graphs.
\end{abstract}


\section{Introduction}\label{intro}

Birnbaum and Sirken (1965) study the situation where patients are sampled indirectly via the hospitals from which they receive treatment. Insofar as a patient may be treated at more than one hospital, the patients are not nested in the hospitals like elements in clustered sampling. Birnbaum and Sirken consider three estimators for such \emph{indirect sampling}. The first one is the classic Horvitz-Thompson (HT) estimator (Horvitz and Thompson, 1952) based on all the sample patients, each of which is weighted by the inverse of the probability of being included in the sample. The second estimator is based on all the sample hospitals and a constructed value for each of them, and the third one is only based on a sub-sample of hospitals determined by a \emph{priority rule}. In particular, the estimator using all the sample hospitals is often referred to as a Hansen-Hurwitz (HH) type estimator. The HH-type estimator and its variations are used for \emph{network sampling} (Sirken, 1970; 2005); it is recast as a ``generalised weight share method'' (Lavalle\`{e}, 2007);  and a modified HH-type estimator is considered for \emph{adaptive cluster sampling} (Thompson, 1990; 1991).  

Zhang and Patone (2017) formally define sampling from finite graphs, in analogy to sampling from finite populations (Neyman, 1934), extending the previous works by Frank (1971, 1980a, 1980b, 2011), which deal with different graph motifs separately. In particular, they  show that each of the aforementioned unconventional sampling techniques can be given different graph sampling representations. Zhang and Oguz-Alper (2020) identify sufficient and necessary conditions for feasible representation of sampling from arbitrary graphs as \textit{bipartite incidence graph sampling (BIGS)}, including indirect, network and adaptive cluster sampling. For instance, the nodes can be the hospitals and the patients and an edge exists between a hospital and any patient that receive treatment at the hospital. This is a bipartite graph since the nodes of the graph are bi-partitioned, where an edge can exist only between two nodes in different parts, but not between any two nodes in the same part.

There are at least three reasons why BIGS is a useful representation of the so-called unconventional sampling techniques in the literature.

First, unconventional sampling techniques are often characterised by the presence of some rules of observation, \emph{in addition} to the probability design of an \emph{initial sample}. For example, under network sampling (Sirken, 1970), ``siblings report each other'' are needed to reach a ``network'' of siblings following an initial sample of households. Under adaptive cluster sampling (Thompson, 1990), sample propagation depends on the ``network'' relationship among the units and the values of the surveyed units. Under graph sampling (Zhang and Patone, 2017), one needs to specify an \emph{observation procedure}, by which the edges of the sample graph are observed following an initial sample of nodes. As demonstrated by Zhang and Oguz-Alper (2020), BIGS can provide a unified representation of various situations of sampling, which are originally described in other terms, where one part of the nodes refer to the initial sampling units and the other part the measurement units of interest, to be referred to as \emph{motifs}, such that the edges represent the observational links between sampling units and motifs. More examples will be given later in this paper.

Next, in the example of indirect sampling of patients via hospitals, one needs to identify all the relevant hospitals outside the initial sample, in order to compute the inclusion probability of a sample patient, which is the information on ``multiplicity'' of sources that must be collected \emph{in addition to} the sample of hospitals and patients. The same requirement exists as well for any other unconventional sampling, such as ``counting rules'' of links between population elements and selection units under network sampling (Sirken, 2005), or the relationship between edge units and their neighbouring networks under adaptive cluster sampling (Thompson, 1990). Under graph sampling, the observation procedure needs to be \emph{ancestral} (Zhang and Patone, 2017), so that one knows which other out-of-sample nodes could have led to the motifs in the sample graph, had they been selected in the initial sample of nodes. The information of multiplicity or ancestry is apparent under BIGS, which is simply the knowledge of the nodes (representing sampling units) that are adjacent to the node representing a sample motif in the BIG.

Finally, as the aim of this paper, one can study design-based estimation under the general setting of \emph{ancestral BIGS} (satisfying the requirement of ancestral observation), where the results are immediately applicable to all the relevant situations. We shall develop a large class of unbiased incidence weighting estimators, based on the sample edges that link the sampling units to the observed motifs. As will be explained, all the three estimators used by Birnbaum and Sirken (1965) are special cases of this class of estimators, which is an insight hitherto unknown in the literature. Many other unbiased estimators can be devised as members of the proposed class, and one can apply the Rao-Blackwell method (Rao, 1945; Blackwell, 1947) to the non-HT estimators, to generate distinct unbiased estimators that can improve the estimation efficiency. Thus, the discovery of the class of incidence weighting estimators provides a potential for efficiency gains.

Below, in Section \ref{basics}, we formally introduce ancestral BIGS, and develop the incidence weighting estimators. The general condition of unbiased estimation is established. New understandings of the three aforementioned estimators are discussed. We consider also the application of Rao-Blackwell method, which motivates a new subclass of the HH-type estimators. Illustrations are given in Section \ref{illustration} of adaptive cluster sampling (Thompson, 1990), line-intercept sampling (Becker, 1991) and simulated graphs, which demonstrate the scope and flexibility of the proposed approach across a variety of situations. Some concluding remarks are given in Section \ref{remark}.

\section{Incidence weighting estimator under BIGS}\label{basics}

Denote by $\mathcal{B} = (F, \Omega; H)$ a bipartite simple directed graph, where $(F, \Omega)$ form a bipartition of the node set $F \cup \Omega$, and each edge in $H$ points from one node in $F$ to another in $\Omega$. No edge exists between any two nodes in $F$ or any two in $\Omega$. For BIGS from $\mathcal{B}$, let $F$ be the set of initial sampling units, and $\Omega$ the population of motifs that are of interest; an edge $(i\kappa)$ that is incident to $i\in F$ and $\kappa \in \Omega$ exists, if and only if the selection of $i$ in a sample $s$ from $F$ leads to the observation of motif $\kappa$ in $\Omega$. The edge set $H$ is unknown to start with. Let the size of $F$ be $M = |F|$, and that of $\Omega$ be $N = |\Omega|$, where $N$ is generally unknown. The incidence relationships corresponding to the edges in $H$ represent thus the observational links between the sampling units and the motifs of interest.

Zhang and Oguz-Alper (2020, Theorem 1) establish the sufficient and necessary conditions, by which an arbitrary instance of graph sampling can be given a feasible BIGS representation. They examine and discuss the BIGS representation of indirect, network and adaptive cluster sampling. For instance, for indirect sampling of patients via hospitals, let $F$ consist of the hospitals and $\Omega$ the patients, where $(i\kappa)\in H$ iff patient $\kappa$ receives treatment at hospital $i$. For network sampling of siblings via households, one can let $F$ consist of the households and $\Omega$ the networks of siblings, i.e. each $\kappa$ represents a group of people who are siblings of each other, where $(i\kappa)\in H$ iff at least one of the siblings in $\kappa$ belongs to household $i$. Adaptive cluster sampling will be discussed in Section \ref{illustration}.

Let $\alpha_i = \{ \kappa : \kappa\in \Omega, (i\kappa) \in H \}$ be the successors of $i$ in $\mathcal{B}$. Given the initial sample $s$ from $F$, the observation procedure of BIGS is \emph{incident} (Zhang and Patone, 2017), such that all the nodes in $\alpha_i$ are included in the sample graph provided $i\in s$; hence, the term BIGS. Let $\Omega_s = \bigcup_{i\in s} \alpha_i$, which consists of all the sample motifs. Following the general definition of sample graph (Zhang and Patone, 2017), the \emph{sample BIG} is given by
\[
\mathcal{B}_s = \big( s, \Omega_s ; H_s\big) \qquad\text{where}\quad H_s = (s\times \Omega_s) \cap H
\]
is the sample of edges. To be able to calculate the inclusion probabilities of each $\kappa$ in $s$, the observation procedure needs to be \emph{ancestral} as well. Let $\beta_{\kappa} = \{ i : i\in F, (i\kappa) \in H \}$ be the \emph{ancestors} (or predecessors) of $\kappa$ in $\mathcal{B}$. Let $\beta(\Omega_s) =  \bigcup_{\kappa \in \Omega_s} \beta_{\kappa}$. The knowledge of ancestry (or multiplicity) amounts to the observation of $\beta(\Omega_s) \setminus s$, although these nodes are not part of the sample graph $\mathcal{B}_s$, such as the out-of-sample hospitals of the sample patients.

\begin{example} \label{ex1} Consider ancestral BIGS from the population BIG below.
\begin{center}
\begin{tikzcd}[cramped]
\kappa_1 & \kappa_2 & \kappa_3 \\
i_1 \arrow[u] & i_2 \arrow[ul] \arrow[u] & i_3 \arrow[u] & i_4
\end{tikzcd}
\end{center}
We have $F = \{ i_1, i_2, i_3, i_4\}$ and $\Omega = \{ \kappa_1, \kappa_2, \kappa_3\}$ and $H = \{ (i_1 \kappa_1), (i_2 \kappa_1), (i_2 \kappa_2), (i_3 \kappa_3) \}$. Suppose $s = \{ i_1, i_3\} \subset F$. By incident observation procedure, we have $\Omega_s = \{ \kappa_1, \kappa_3\}$ and $H_s = \{ (i_1 \kappa_1), (i_3\kappa_3)\}$, and the sample graph $\mathcal{B}_s = (s, \Omega_s; H_s)$ as defined above. In addition, we observe $\beta(\Omega_s) \setminus s = \{ i_2\}$, where $i_2$ is not part of the sample BIG.
\end{example}

\subsection{The incidence weighting estimator}

Let $y_{\kappa}$ be an unknown constant associated with motif $\kappa$, for $\kappa\in \Omega$, given the population graph $\mathcal{B}$. The aim is to estimate the total $\theta = \sum_{\kappa\in \Omega} y_{\kappa}$, including e.g. $y_{\kappa} \equiv 1$. Given the sample graph $\mathcal{B}_s$, let $\{ W_{i\kappa}; (i\kappa)\in H_s\}$ be the \emph{incidence weights} of the sample edges, and $W_{i\kappa}\equiv 0$ if $(i\kappa)\not \in H_s$. The \emph{incidence weighting estimator (IWE)} is given by
\begin{equation} \label{IWE}
\hat{\theta} = \sum_{(i\kappa) \in H_s} W_{i\kappa} y_{\kappa}/\pi_{(i\kappa)}
\end{equation}
where $\pi_{(i\kappa)} = \mbox{Pr}\big( (i\kappa) \in H_s\big)$. Under BIGS, we have $\pi_{(i\kappa)} = \pi_i = \mbox{Pr}(i \in s) = E(\delta_i)$, where $\delta_i = 1$ or 0 indicates if $i\in s$ or not, and $\pi_i$ is the probability of $\delta_i = 1$, for $i\in F$. Notice that the definition \eqref{IWE} allows for sample dependent weights $W_{i\kappa}$.

\paragraph{\em Proposition 1.} The IWE by \eqref{IWE} is unbiased for $\theta$ provided, for each $\kappa\in \Omega$,
\begin{equation}\label{unbias}
 \sum_{i \in \beta_{\kappa}} E(W_{i\kappa} | \delta_i = 1) = 1 ~.
\end{equation}

\begin{proof}
The expectation of $\hat{\theta}$ with respect to the sampling distribution of $s$ is given by
\begin{align*}
E(\hat{\theta}) & = \sum_{i \in F} \frac{E(\delta_i)}{\pi_i} E(\sum_{\kappa\in \alpha_i} W_{i\kappa} y_{\kappa} | \delta_i =1)
%= \sum_{i \in F} \sum_{\kappa \in \alpha_i} E(W_{i\kappa} | \delta_i = 1) y_{\kappa} \\
= \sum_{\kappa\in \Omega} y_{\kappa}  \sum_{i \in \beta_{\kappa}} E(W_{i\kappa} | \delta_i = 1) = \theta
\end{align*}
since $\pi_{(i\kappa)} = \pi_i$ under BIGS, and $\sum_{i \in \beta_{\kappa}} E(W_{i\kappa} | \delta_i = 1) =1$ by stipulation.
\end{proof}

The condition \eqref{unbias} ensures that the IWE is unbiased under repeated sampling. When the weights are constant of sampling, denoted by $\omega_{i\kappa}$ for distinction, it reduces to $\sum_{i\in \beta_{\kappa}} \omega_{i\kappa} = 1$ for any $\kappa \in \Omega$. Let $\pi_{ij}$ be the second-order sample inclusion probability of $i,j\in F$.

\paragraph{\em Proposition 2.} The BIG sampling variance of an unbiased IWE is given by
\begin{equation}\label{vIWE}
V(\hat{\theta}) = \sum_{\kappa \in \Omega}\sum_{\ell \in \Omega} (\Delta_{\kappa \ell} -1) y_{\kappa} y_{\ell}
\end{equation}
where
\[
\Delta_{\kappa \ell} = \sum_{i \in \beta_{\kappa}} \sum_{j \in \beta_{\ell}} \frac{\pi_{ij}}{\pi_i \pi_j}
E\big( W_{i\kappa}W_{j\ell} | \delta_i \delta_j =1\big) ~.
\]

\begin{proof} Given unbiased $\hat{\theta}$, we have $V(\hat{\theta}) = E(\hat{\theta}^2) - \theta^2$, where
\begin{align*}
E\big( \hat{\theta}^2 \big) &= \sum_{\kappa \in \Omega} \sum_{\ell \in \Omega} y_{\kappa} y_{\ell} 
\sum_{i\in \beta_{\kappa}} \sum_{j\in \beta_{\ell}} E\Big( \frac{\delta_i \delta_j}{\pi_i \pi_j} W_{i\kappa} W_{j\ell} \Big) \\
%&= \sum_{\kappa \in \Omega} \sum_{\ell \in \Omega} y_{\kappa} y_{\ell} 
%\sum_{i\in \beta_{\kappa}} \sum_{j\in \beta_{\ell}} \frac{E(\delta_i \delta_j =1)}{\pi_i \pi_j} E(W_{i\kappa} W_{j\ell}  | \delta_i \delta_j =1) \\
&= \sum_{\kappa \in \Omega} \sum_{\ell \in \Omega} y_{\kappa} y_{\ell} 
\sum_{i\in \beta_{\kappa}} \sum_{j\in \beta_{\ell}}  \frac{\pi_{ij}}{\pi_i \pi_j} E(W_{i\kappa} W_{j\ell}  | \delta_i \delta_j =1)
\end{align*}
since $W_{i\kappa} W_{j\ell} =0$ if $\delta_i \delta_j =0$ under BIGS, for any $(i\kappa), (j\ell) \in H$. The result follows now from taking the difference of $E(\hat{\theta}^2)$ and $\theta^2 = (\sum_{\kappa\in \Omega} y_{\kappa})^2$.
\end{proof}

\subsection{HT-type estimator}

Let $\pi_{(\kappa)} = \mbox{Pr}(\kappa \in \Omega_s)$ and $\pi_{(\kappa\ell)} = \mbox{Pr}(\kappa \in \Omega_s, \ell \in \Omega_s)$ for $\kappa, \ell\in \Omega$, where parentheses are used in the subscript to distinguish these inclusion probabilities of the motifs from those of the sampling units. The HT-estimator is given by
\[
\hat{\theta}_y = \sum_{\kappa \in \Omega_s} y_{\kappa}/\pi_{(\kappa)}
\]
where $V(\hat{\theta}_y) = \sum_{\kappa \in \Omega} \sum_{\ell \in \Omega} (\pi_{(\kappa\ell)} / \pi_{(\kappa)} \pi_{(\ell)} -1) y_{\kappa} y_{\ell}$. Under BIGS, we have
\begin{gather*}
\pi_{(\kappa)} = 1 - \bar{\pi}_{\beta_{\kappa}} = 1 - \mbox{Pr}\big( \beta_{\kappa} \cap s = \emptyset \big) \\
\pi_{(\kappa \ell)} = 1 - \big( \bar{\pi}_{\beta_{\kappa}} + \bar{\pi}_{\beta_{\ell}} - \bar{\pi}_{\beta_{\kappa} \cup \beta_{\ell}} \big) ~.
\end{gather*}
where $\bar{\pi}_{\beta_k}$ is the exclusion probability of $\beta_k$ in $s$, which is the probability that none of the ancestors of $\kappa$ in $\mathcal{B}$ is included in the initial sample $s$, and the knowledge of the out-of-sample ancestors $\beta_{\kappa} \setminus s$ is required to compute $\bar{\pi}_{\beta_{\kappa}}$. Similarly for $\bar{\pi}_{\beta_{\kappa} \cup \beta_{\ell}}$.

The HT-estimator is a special case of the IWE, where the weights $W_{i\kappa}$ satisfy
\begin{equation}\label{Yhat}
\sum_{i \in s \cap \beta_{\kappa}} W_{i\kappa}/\pi_i = 1/\pi_{(\kappa)} ~.
\end{equation}
Notice that these weights $W_{i\kappa}$ are not constant of sampling if $|\beta_{\kappa}| > 1$, since they depend on how $s$ intersects $\beta_{\kappa}$. Take $\kappa_1$ in Example before, where $\beta_{\kappa_1} = \{ i_1, i_2\}$. By \eqref{Yhat}, we have $W_{i_1 \kappa_1} = \pi_{i_1} / \pi_{(\kappa_1)}$ if $s\cap \{ i_1, i_2\} = \{ i_1\}$, and $W_{i_2 \kappa_1} = \pi_{i_2} / \pi_{(\kappa_1)}$ if $s\cap \{ i_1, i_2\} = \{ i_2\}$, and one can e.g. let $(W_{i_1\kappa},  W_{i_2 \kappa_1}) = \frac{1}{2} (\pi_{i_1} , \pi_{i_2})/\pi_{(\kappa_1)}$ if $s\cap \{ i_1, i_2\} = \{ i_1, i_2\}$.

To see that the weights given by \eqref{Yhat} satisfy the condition \eqref{unbias} generally, let $\phi_{s_{\kappa}}$ be the probability that the \emph{sample intersection} is $s_{\kappa} = s\cap \beta_{\kappa}$ for $\kappa\in \Omega$, where $\pi_{(\kappa)} = \sum_{s_{\kappa}} \phi_{s_{\kappa}}$ over all possible $s_{\kappa}$. Given \eqref{Yhat}, for any $\kappa \in \Omega$, we have then
\[
\sum_{i \in \beta_{\kappa}} E(W_{i\kappa} | \delta_i =1)
= \sum_{i \in \beta_{\kappa}} \sum_{s_{\kappa} \ni i} \frac{\phi_{s_{\kappa}}}{\pi_i} W_{i\kappa}
= \sum_{s_{\kappa}} \phi_{s_{\kappa}} \sum_{i \in s_{\kappa}} \frac{W_{i\kappa}}{\pi_i}
= \sum_{s_{\kappa}} \frac{\phi_{s_{\kappa}}}{\pi_{(\kappa)}} = 1~.
\]
Arguing similarly in terms of the joint probability that the sample intersections for $\kappa$ and $\ell$ are $s_{\kappa}$ and $s_{\ell}$, it can be shown that $\Delta_{\kappa \ell}$ in \eqref{vIWE} reduces to $\pi_{(\kappa\ell)} / \pi_{(\kappa)} \pi_{(\ell)}$ given \eqref{Yhat} and \eqref{unbias}.

More generally, let $\eta_{s_{\kappa}} = \pi_{(\kappa)} \sum_{i\in s_{\kappa}} W_{i\kappa}/\pi_i$ for any weights $W_{i\kappa}$ that are not constants of sampling. To satisfy the condition \eqref{unbias}, for any $\kappa\in \Omega$, the weights must be such that
\begin{equation} \label{W}
\sum_{s_{\kappa}} \phi_{s_{\kappa}} \eta_{s_{\kappa}} = \pi_{(\kappa)} ~.
\end{equation}
The HT-estimator is the special case where $\eta_{s_{\kappa}} \equiv 1$. It is possible to assign $\eta_{s_{\kappa}}$ that differs from 1 for different sample intersects $s_{\kappa}$, subjected to the restriction \eqref{W}. Any estimator satisfying \eqref{W} but not \eqref{Yhat} may be referred to as a \emph{HT-type estimator}.

\subsection{HH-type estimator}

While a HT-type estimator uses sample dependent weights $W_{i\kappa}$,
a HH-type estimator uses weights $\omega_{i\kappa}$ that are constant of sampling. The condition \eqref{unbias} is reduced to $\sum_{i\in \beta_{\kappa}} \omega_{i\kappa} = 1$, for any $\kappa \in \Omega$. Birnbaum and Sirken (1965) observe that 
\[
\sum_{\kappa \in \Omega} y_{\kappa} = \sum_{\kappa \in \Omega} y_{\kappa} \sum_{i\in \beta_{\kappa}} \omega_{i\kappa} =
\sum_{i\in F} \sum_{\kappa \in \alpha_i} \omega_{i\kappa} y_{\kappa} ~.
\]
It follows that the HH-type estimator given by
\begin{equation}\label{Zhat}
\hat{\theta}_z = \sum_{i \in s} z_i / \pi_i \qquad\text{and}\qquad z_i = \sum_{\kappa \in \alpha_i} \omega_{i\kappa} y_{\kappa}
\end{equation}
is unbiased for $\theta$ under repeated sampling, where $z_i$ is a constructed constant for each initial sample unit $i$. The BIG sampling variance of $\hat{\theta}_z$ is given by
\[
V(\hat{\theta}) = \sum_{i \in F} \sum_{j \in F} \big( \frac{\pi_{ij}}{\pi_i \pi_j} -1) z_i z_j ~.
\]

Notice that one only needs $z_i$ for the initial sample units in order to apply $\hat{\theta}_z$, which is possible provided ancestral BIGS. Moreover, the HH-type estimator \eqref{Zhat} defines actually a family of estimators, depending on the choice of $\omega_{i\kappa}$, although Birnbaum and Sirken (1965) use only the equal weights $\omega_{i\kappa} = 1/|\beta_{\kappa}|$. The corresponding $\hat{\theta}_z$ is referred to as the \emph{multiplicity estimator}, denoted by $\hat{\theta}_{z\beta}$. Variations of the multiplicity estimator under other settings of indirect, network sampling are considered by Sirken (1970), Sirken and Levy (1974), Sirken (2004) and Lavalle\`{e} (2007). Unlike the HT-estimator, it is in principle possible to apply the Rao-Blackwell method to improve the HH-type estimator, to which we return in Section \ref{RB}. Some other HH-type estimators will be discussed then.

\subsection{Priority-rule estimator}

Birnbaum and Sirken (1965) invent a third estimator based on a \emph{prioritised subset} of $H_s$, where they let $I_{i\kappa} =1$ if $i = \min\big( s \cap \beta_{\kappa} \big)$ and 0 otherwise, i.e. if unit $i$ happens to be enumerated first in the frame $F$ among all the in-sample ancestors of $\kappa$, for each $\kappa\in \Omega_s$. Their \emph{priority-rule estimator} based on $\{ (i\kappa) : I_{i\kappa} = 1, (i\kappa)\in H_s\}$ is given by
\begin{equation} \label{Phat}
\hat{\theta}_p = \sum_{(i\kappa) \in H_s} \frac{I_{i\kappa} \omega_{i\kappa} y_{\kappa}}{p_{i\kappa} \pi_i}
\end{equation}
where $p_{i\kappa} =\mbox{Pr}\big( I_{i\kappa} =1 | (i\kappa) \in H_s \big) =\mbox{Pr}\big( I_{i\kappa} =1 | \delta_i =1\big)$ is the conditional probability that $(i\kappa)$ is prioritised given $(i\kappa) \in H_s$, and $\omega_{i\kappa} = 1/|\beta_{\kappa}|$ are the equal weights for any $\kappa\in \Omega$. Clearly, other priority rules or choices of $\omega_{i\kappa}$ are possible.

One can easily recognise $\hat{\theta}_p$ as a special case of IWE with $W_{i\kappa} =  I_{i\kappa} \omega_{i\kappa}/ p_{i\kappa}$. It can satisfy the unbiasedness condition \eqref{unbias},  \emph{provided} $p_{i\kappa} >0$ for all $(i\kappa) \in H_s$, in which case $E(W_{i\kappa} | \delta_i =1) = \omega_{i\kappa}$. Birnbaum and Sirken (1965) did not provide an expression of $V(\hat{\theta}_p)$, but indicated that it is unwieldy. Now that $\hat{\theta}_p$ is a special case of IWE, its variance follows readily from Proposition 2. Let $p_{i\kappa, j\ell} = \mbox{Pr}(I_{i\kappa} I_{j\ell} =1 | \delta_i \delta_j =1)$, we have
\[
\Delta_{\kappa \ell} = \sum_{i \in \beta_{\kappa}} \sum_{j \in \beta_{\ell}} \frac{\pi_{ij} p_{i\kappa, j\ell}}{\pi_i \pi_j p_{i\kappa} p_{j\ell}}
\omega_{i\kappa} \omega_{j\ell}
\]
in \eqref{vIWE}, such that
\[
% = \sum_{\kappa \in \Omega}\sum_{\ell \in \Omega} \Big( \sum_{i \in \beta_{\kappa}} \sum_{j \in \beta_{\ell}} \frac{\pi_{ij} p_{i\kappa, j\ell}}{\pi_i \pi_j p_{i\kappa} p_{j\ell}}  -1 \Big) y_{\kappa} y_{\ell} \\
V(\hat{\theta}_p) = \sum_{(i\kappa)\in H} \sum_{(j\ell)\in H} \Big( \frac{\pi_{ij} p_{i\kappa, j\ell}}{\pi_i \pi_j p_{i\kappa} p_{j\ell}} - 1 \Big)
\omega_{i\kappa} \omega_{j\ell} y_{\kappa} y_{\ell}
\]
because $\sum_{i\in \beta_{\kappa}} \omega_{i\kappa} =1$ for any $\kappa \in \Omega$. An unbiased variance estimator can be given by
\[
\widehat{V}(\hat{\theta}_p) = \sum_{(i\kappa)\in H_s} \sum_{(j\ell)\in H_s}  \Big( \frac{\pi_{ij} p_{i\kappa, j\ell}}{\pi_i \pi_j p_{i\kappa} p_{j\ell}} - 1 \Big) \frac{\omega_{i\kappa} \omega_{j\ell}}{\pi_{ij}} y_{\kappa} y_{\ell} ~.
\]
The priority probabilities $p_{i\kappa}$ and $p_{i\kappa, j\ell}$ depend on the priority rule, as well as the sampling design. The details for the estimtator of Birnbaum and Sirken (1965) under initial simple random sampling (SRS) without replacement of $s$ are given in Appendix \ref{pprob}.

It should be noticed that the priority rule is not part of sampling; the sample graph $\mathcal{B}_s$ includes all the edges incident to every sample unit in $s$. Had one applied subsampling by randomly selecting one of the edges incident to each $i$ in $s$ with some designed probabilities, the sample graph would have contained one and only one edge from each sample unit. Instead, the priority rule selects only one sample edge incident to each motif in $\Omega_s$ for the purpose of estimation. There is a possibility that a unit $i$ can be sampled but never prioritised, in which case $\hat{\theta}_p$ would be biased. For an extreme example, suppose a motif $\kappa$ is incident to all the sampling units in $F$, then the last unit in $F$ can never be prioritised (for $\kappa$) according to the priority rule of Birnbaum and Sirken (1965). Generally, $\hat{\theta}_p$ is \emph{biased} under this priority rule, provided there exists at least one motif $\kappa$ in $\Omega$, where
\[
|\beta_{\kappa}|>1 \quad\text{and}\quad \mbox{Pr}( |s_{\kappa}| \leq 1) =0
\]
such that the ancestor $i = \max(\beta_{\kappa})$ has no chance of being prioritised when it is in $s$. The probability above depends on the ordering of sampling units in $F$, as well as the initial sample size. Given any ordering of the units in $F$, as the initial sample increases, it is possible for $\hat{\theta}_p$ to behave more erratically and become biased eventually.

\subsection{Using Rao-Blackwell method} \label{RB}

The minimal sufficient statistic under BIGS is $\{ (\kappa, y_{\kappa}) : \kappa \in \Omega_s \}$, or simply $\Omega_s$ as long as one keeps in mind that the $y$-values are constants associated with the motifs. Let $\hat{\theta}$ be an unbiased IWE. Applying the Rao-Blackwell method to $\hat{\theta}$ yields $\hat{\theta}_{RB} = E(\hat{\theta} | \Omega_s)$ as an improved estimator, if the conditional variance $V(\hat{\theta} | \Omega_s)$ is positive. Since the HT-estimator $\hat{\theta}_y$ is fixed conditional on $\Omega_s$, we have $\hat{\theta}_{yRB} \equiv \hat{\theta}_y$. For a non-HT estimator, it is in principle possible that the RB method can improve its efficiency, as illustrated below.

\begin{example}[cont'd] Given $|s| =1$, there are 4 distinct initial samples, leading to 4 distinct $\Omega_s$ under BIGS, such that  $V(\hat{\theta} | \Omega_s) = 0$ and $\hat{\theta}_{RB} = \hat{\theta}$ for any unbiased IWE. Given $|s| = 2$, there are 6 different initial samples, leading to 5 distinct $\Omega_s$, where both $s = \{ i_1, i_2\}$ and $s' = \{ i_2, i_4\}$ lead to the same motifs $\{ \kappa_1, \kappa_2\}$, so that $\hat{\theta}_{RB} \neq \hat{\theta}$ given motif sample $\{ \kappa_1, \kappa_2\}$, if $\hat{\theta}(s) \neq \hat{\theta}(s')$. Take e.g. the HH-type estimator $\hat{\theta}_z$ by \eqref{Zhat}, we have
\begin{gather*}
\hat{\theta}_z(s) = \frac{\omega_{i_1 \kappa_1}}{\pi_{i_1}} y_{\kappa_1} + \frac{\omega_{i_2 \kappa_1}}{\pi_{i_2}} y_{\kappa_1}
+  \frac{\omega_{i_2 \kappa_2}}{\pi_{i_2}} y_{\kappa_2} \neq 
\hat{\theta}_z(s') = \frac{\omega_{i_2 \kappa_1}}{\pi_{i_2}} y_{\kappa_1} + \frac{\omega_{i_2 \kappa_2}}{\pi_{i_2}} y_{\kappa_2} \\
\hat{\theta}_{zRB} = \frac{p(s)}{p(s) + p(s')} \cdot \frac{\omega_{i_1 \kappa_1}}{\pi_{i_1}} y_{\kappa_1}
+ \frac{\omega_{i_2 \kappa_1}}{\pi_{i_2}} y_{\kappa_1} +  \frac{\omega_{i_2 \kappa_2}}{\pi_{i_2}} y_{\kappa_2} ~.
\end{gather*}
\end{example}

The calculation required for the RB method may be intractable, if the conditional sample space of $s$ given $\Omega_s$ is large and the initial sampling design $p(s)$ is not fully specified, which is common in practice for designs with unequal inclusion probabilities over $F$. Moreover, the result of RB method is generally not a unique minimum variance unbiased estimator under BIGS, because the minimal sufficient statistic is not complete. It is thus worth exploring other useful choices of the IWE. Due to the inherent shortcoming of the priority-rule estimator pointed out earlier, we concentrate on the HH-type estimator $\hat{\theta}_z$ below.

Consider the special case where $|\alpha_i| \equiv 1$ in the population BIG, such as when sampling households via persons. Suppose first with-replacement sampling of $s$, where the different draws generate an IID sample, and compare $\hat{\theta}_y$ and $\hat{\theta}_z$ based on a single draw. Let $p_i$ and $p_{(\kappa)} = \sum_{i\in \beta_{\kappa}} \pi_i$ be the respective \emph{selection} probabilities. We have $p_{ij} = p_i$ if $i=j$ and $0$ if $i\neq j$, and $p_{(\kappa \ell)} = p_{(\kappa)}$ if $\kappa=\ell$ and $0$ if otherwise, now that $|\alpha_i| \equiv 1$. We have
\[
V(\sum_{i\in s} z_i/p_i) - V(\sum_{\kappa\in \Omega_s} y_{\kappa}/p_{(\kappa)}) =
\sum_{\kappa\in \Omega} \Big( \sum_{i\in \beta_{\kappa}} \omega_{i\kappa}^2/p_i - 1/p_{(\kappa)} \Big) y_{\kappa}^2 =0
\]
if $\omega_{i\kappa} \equiv p_i/p_{(\kappa)}$ for $i\in \beta_{\kappa}$, given which we have $\hat{\theta}_z = \hat{\theta}_{zRB}$. The variance of any other $\hat{\theta}_z$ would be larger, as long as $\omega_{i\kappa}/p_i$ is not a constant over $\beta_{\kappa}$, because
\[
E\big( \frac{z_i}{p_i} | \kappa\in \Omega_s \big) = \sum_{i\in \beta_{\kappa}} \frac{p_i w_{i\kappa} y_{\kappa}}{p_{(\kappa)} p_i} = \frac{y_{\kappa}}{p_{(\kappa)}} \quad\text{and}\quad 
V\big( \frac{z_i}{p_i} | \kappa\in \Omega_s\big) = y_{\kappa}^2 V\big( \frac{\omega_{i\kappa}}{p_i} | \kappa\in \Omega_s \big) > 0 ~.
\]

A similar argument holds approximately for the choice $\omega_{i\kappa} \propto \pi_i$ under sampling without replacement of $s$, provided $\pi_{ij} \approx \pi_i \pi_j$ and $\pi_{(\kappa \ell)} \approx  \pi_{(\kappa)} \pi_{(\ell)}$, as in the case of sampling households via persons with a small sampling fraction $|s|/|F|$. This can make $z_i/\pi_i$ more similar to each other over $F$, which is advantageous with respect to the anticipated mean squared error of $\hat{\theta}_z$ under the sampling design and a population model of $z_i$, according to Theorem 6.2 of Godambe and Joshi (1965).

To make $z_i/\pi_i$ more similar to each other over $F$ without the restriction $|\alpha_i | =1$, one may consider setting $\omega_{i\kappa} < \omega_{j\kappa}$ if $|\alpha_i| > |\alpha_j|$, despite $\pi_i = \pi_j$, because there are more motifs contributing to $z_i$ than $z_j$. Thus, under general unequal-probability sampling of $s$, it may be reasonable to consider the \emph{probability and inverse-degree adjusted  (PIDA)} weights
\begin{equation} \label{paid}
\omega_{i\kappa} \propto \pi_i/|\alpha_i|^{\gamma}
\end{equation}
subjected to the condition \eqref{unbias}, where $\gamma>0$ is a tuning constant of choice. Denote by $\hat{\theta}_{z\alpha\gamma}$ the corresponding PIDA-IWE. The multiplicity estimator $\hat{\theta}_{z\beta}$ becomes a special case of $\hat{\theta}_{z\alpha\gamma}$ given $\gamma =0$ and constant $\pi_i$ over $F$.

Notice that to apply the weights \eqref{paid} with $\gamma\neq 0$, one needs to know $|\alpha_i|$ for all $i\in \beta_{\kappa}$ and $\kappa\in \Omega_s$, in addition to the ancestral observation of $\beta_{\kappa}$. For instance, under indirect sampling of children via parents, one would need to collect the number of children for the out-of-sample parents in $\beta(\Omega_s)\setminus s$ as well. For network sampling of siblings via households, one would need to collect the number of other sibling networks in each household $i$ with at least one member from a sample sibling network $\kappa$.

\section{Illustrations} \label{illustration}

\subsection{Adaptive cluster sampling}

Consider the example of adaptive cluster sampling (ACS) discussed by Thompson (1990). The population $F$ consists of 5 grids, with $y$-values $\{ 1, 0, 2, 10,1000\}$. Each grid has either one or two neighbours which are adjacent in the given sequence, as in the graph $G$ below, where as Thompson (1990) we simply denote each grid by its $y$-value.
\begin{center}
\begin{tikzcd}[cramped]
G: & 1\arrow[r] & 0\arrow[l]\arrow[r] & 2 \arrow[l] \arrow[r]& 10 \arrow[l]\arrow[r] & 1000\arrow[l]
\end{tikzcd}
\end{center}
Given an initial sample of size 2 by SRS from $F$, one would survey all the neighbour grids (in both directions if possible) of a sample grid $i$ if $y_i$ exceeds the threshold value 5 but not otherwise. The observation procedure is repeated for all the neighbour grids, which may or may not generate further grids to be surveyed. The process is terminated, when the last observed grids are all below the threshold. The interest is to estimate the total amount of species (or mean per grid) over the given area.

In particular, the grid 2 is a so-called \emph{edge unit}, which can be observed from 10 or 1000, but would not lead to 10 or 1000 if only 2 is selected in $s$. The inclusion probability of grid $2$ under ACS cannot be calculated correctly when it is selected in $s$ but not 10 or 1000, in which case the knowledge of multiplicity (or ancestry) is lacking. Thompson (1990) proposes a modified HT-estimator which uses the grid 2 in estimation, only if it is selected on its own, the probability of which is known from the design of the initial sample.

Zhang and Oguz-Alper (2020) develop feasible BIGS representations of ACS from $G$ above. Here we use one of them to illustrate how the IWE can be applied to ACS. The population BIG is given by $\mathcal{B} = (F, F; H)$, with $\Omega = F$ and $H$ as below.
\begin{center}
\begin{tikzcd}[cramped]
\mathcal{B}: & 1 & 0 & 2 & 10 & 1000 \\
& 1 \arrow[u] & 0 \arrow[u] & 2 \arrow[u] & 10 \arrow[u]\arrow[ur] & 1000 \arrow[ul]\arrow[u]
\end{tikzcd}
\end{center}
Zhang and Oguz-Alper (2020) point out that it is possible to consider BIGS from $\mathcal{B}$, where the observational links between $(10, 2)$ and $(1000, 2)$ under ACS are removed to ensure ancestral observation, and apply the classic HT-estimator under this BIGS representation of ACS from $G$. They show that the two strategies (ACS, modified HT) and (BIGS, HT) actually lead to the same estimator. The difference is that one cannot apply the RB method to the HT-estimator under BIGS, as one can with the modified HT-estimator under ACS. We refer to Zhang and Oguz-Alper (2020) for more details.

Thompson (1990) proposes also a modified HH-type estimator, where an edge unit is used in estimation only if it is selected in $s$ directly. This modified HH-type estimator is simply the multiplicity estimator $\hat{\theta}_{z\beta}$ under BIGS from $\mathcal{B}$, with equal weights $\omega_{i\kappa} = 1/|\beta_{\kappa}|$ in \eqref{Zhat}. The two strategies (ACS, modified HH-type) and (BIGS, $\hat{\theta}_{z\beta}$) lead to the same estimator. Moreover, application of the RB method to $\hat{\theta}_{z\beta}$ is the same as that for the modified HH-type estimator; we refer to Thompson (1990) for the details.

Finally, since the contiguous grids that form a network are all observed together under ACS if any of them is observed, ancestral BIGS from $\mathcal{B}$ entails the observation of $|\alpha_i|$ needed for the PIDA weights given by \eqref{paid}. However, since $|\alpha_i |$ is the same for all the grids in the same network and the initial sampling is SRS, the weights by \eqref{paid} are all equal in this case, so that the estimator $\hat{\theta}_{z\alpha\gamma}$ coincides with the multiplicity estimator $\hat{\theta}_{z\beta}$.

\subsection{Line-intercept sampling}

Line-intercept sampling (LIS) is a method of sampling habitats in a region, where a habitat is sampled if a chosen line segment transects it. The habitat may e.g. be animal tracks, roads, forestry, which are of irregular shapes. Kaiser (1983) considers the general situation, where a point is randomly selected on the map and an angle is randomly chosen, yielding a line segment of fixed length or transecting the whole area in the chosen direction. Repetition generates an IID sample of lines. In the simplest setting, each transect line is selected at random by selecting randomly a position along a fixed \emph{baseline} that traverses the whole study area, in the direction perpendicular to the baseline. We apply IWE under BIGS to the following example of LIS (Becker, 1991) under this simple setting.

\begin{figure}[ht]
\resizebox{163mm}{104mm}{\includegraphics{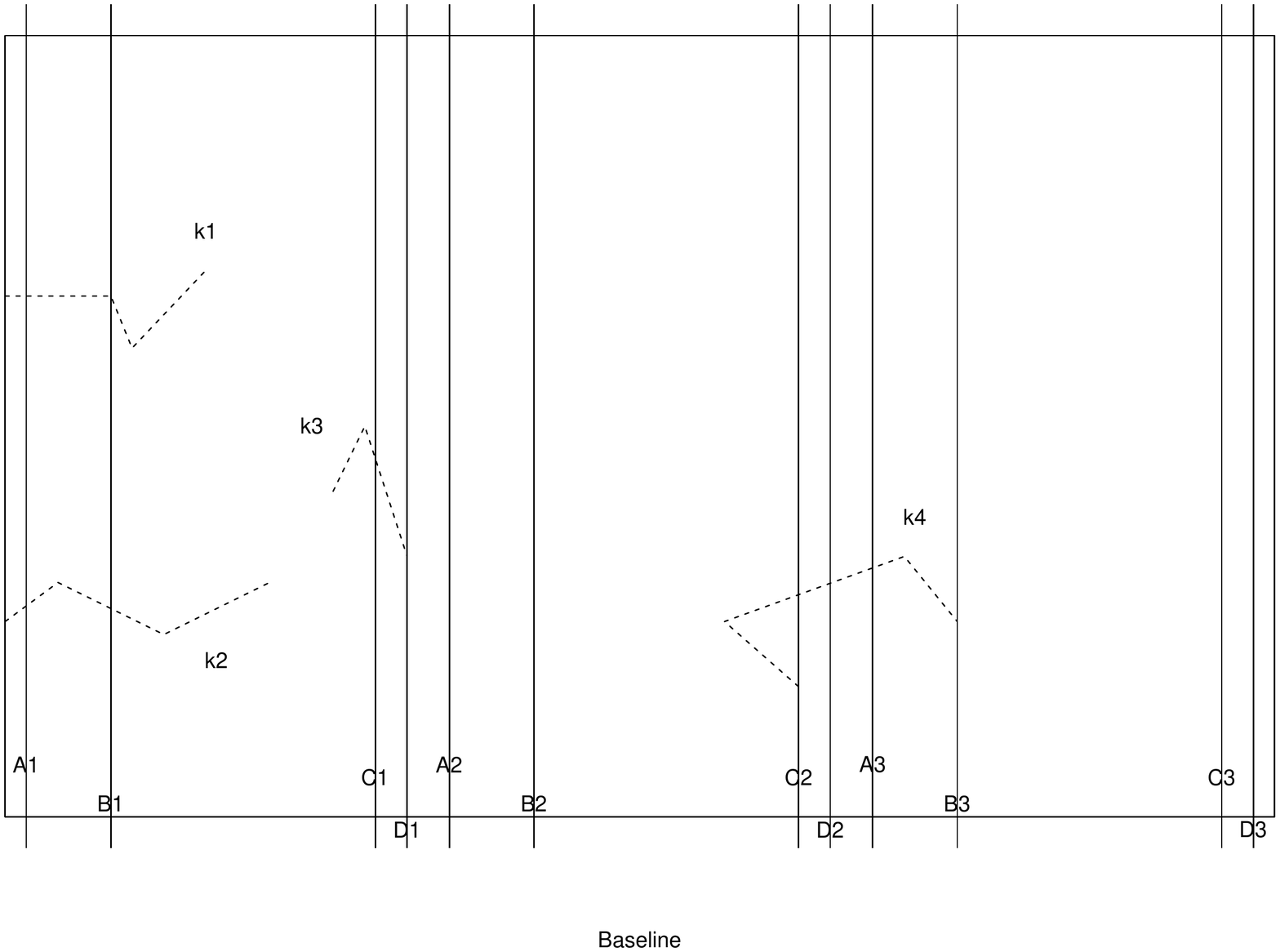}}
\caption{LIS for wolverinne with 4 systematic samples (A, B, C, D) of 3 positions each.}  \label{fig-LIS}
\end{figure} \label{fig-LIS}

The aim is to estimate the total number of wolverines in the mapped area, as sketched in Figure \ref{fig-LIS}. Four systematic samples A, B, C and D, each containing 3 positions, are drawn on the baseline that is equally divided into 3 segments of length 12 miles each. Following the 12 selected lines \emph{and} any wolverine track that intercepts them yields 4 observed tracks, denoted by $\kappa = 1, ..., 4$ and heuristically indicated by the dashed lines in Figure \ref{fig-LIS}. Let $y_{\kappa}$ be the associated number of wolverines, and $L_{\kappa}$ the length of the projection of $\kappa$ on the baseline. From top to bottom and left to right, we observe $(y_1, L_1) = (1, 5.25)$, $(y_2, L_2) = (2, 7.5)$, $(y_3, L_3) = (2, 2.4)$ and $(y_4, L_4) = (1, 7.05)$.

\subsubsection{Feasible BIGS representation of LIS}

First we construct a feasible BIGS representation of LIS in this case. Given the observed tracks, partition the baseline into 7 \emph{projection segments}, each with associated length $x_i$, for $i=1, ..., 7$ from left to right, where $x_1$ refers to the overlapping projection of $\kappa =1$ and 2, $x_2$ the projection of $\kappa =2$ that does not overlap with $\kappa =1$, $x_3$ the distance between projections of $\kappa =2$ and 3, $x_4$ the projection of $\kappa = 3$, $x_5$ the distance between projections of $\kappa =3$ and 4, $x_6$ the projection of $\kappa =4$, and $x_7$ the distance between $\kappa =4$ and right-hand border. The probability that the $i$-th projection segment is selected by a systematic sample is $p_i = x_i/12$. The 4 systematic samples are IID.

The sample BIG on the $r$-th draw is given by $\mathcal{B}_r = (s_r, \Omega_r; H_r)$, where $s_r$ contains the selected projection segments, and $\alpha_i$ the wolverine tracks that intercept the sampled line originating from $i\in s_r$, such that $\Omega_r = \bigcup_{i\in s_r} \alpha_i$ and $H_r = \bigcup_{i\in s_r}  i\times \alpha_i$. In this example, we have $s_1 = s_2 =\{ 1, 5, 6\}$, yielding $\Omega_1= \Omega_2 = \{ 1, 2, 4\}$ on the first two draws A and B, and $s_3 = s_4 = \{ 4, 6, 7\}$, yielding $\Omega_3 = \Omega_4 = \{ 3, 4\}$ on the last two draws C and D. The distinct projection segments selected over all the draws are $s = \bigcup_{r=1}^4 s_r = \{ 1,4,5,6,7\}$, and the distinct tracks are $\Omega_s = \bigcup_{r=1}^4 \Omega_r = \{ 1, 2, 3, 4\}$. Let $F^* = \{ 1, 2, ..., 7\}$ contain the 7 projection segments constructed from $(s, \Omega_s)$, and $H_s = \bigcup_{r=1}^4 H_r$. Let $\mathcal{B}^* = (F^*, \Omega_s; H_s)$ be given as below:
\begin{center}
\begin{tikzcd}[cramped]
\kappa=1 & \kappa=2 & & \kappa=3 & & \kappa=4 \\
i=1 \arrow[u] \arrow[ur] & i=2 \arrow[u] & i=3 & i=4 \arrow[u] & i=5 & i=6 \arrow[u] & i=7
\end{tikzcd}
\end{center}
Let $\beta_{\kappa}^*$ be the ancestors of $\kappa$ in $\mathcal{B}^*$, where $\beta_1^* = \{ 1\}$, $\beta_2^* = \{ 1, 2\}$, $\beta_3^* = \{ 4\}$ and $\beta_4^* = \{ 6\}$.

Let $\Omega = \{ 1, ..., \kappa, ..., N\}$ contain all the wolverine tracks in the area, where $N \geq 4$ given the sample $\Omega_s$. Let $F(\Omega) = \{ 1, ..., i, ..., M\}$ be the sampling frame, which consists of all the projection segments constructed from $\Omega$. Let $H = \{ (i\kappa); i\in F, \kappa\in \Omega\}$, where an edge exists from $i$ to $\kappa$ provided $\kappa$ intercepts any line that originates from the $i$-th projection segment. The population BIG is given by $\mathcal{B} = (F, \Omega; H)$. By Theorem 1 of Zhang and Oguz-Alper (2020) the LIS can be represented as BIGS from $\mathcal{B}$ where, in particular, the observation procedure of LIS ensures that BIGS from $\mathcal{B}$ is ancestral for $\Omega_s$.

Notice that in the population frame $F(\Omega)$ based on $\Omega$, one may need to further partition the projection segments in $F^*$, in order to accommodate the edges and the ancestor sets of the unobserved tracks in $\Omega \setminus \Omega_s$. For instance, suppose there is another track that can only be intercepted from the $7$-th projection segment in $F^*$ and the track does not reach the right-hand border, then this projection segment should be partitioned into 3 segments in $F$ and $(F, H)$ would differ from $(F^*, H_s)$ accordingly. Nevertheless, each observed track $\kappa$ in $\Omega_s$ can only be intercepted from any projection segment $i$ in $F^*$ with $(i\kappa) \in H_s$, such that its ancestor set in $\mathcal{B}$ is observed in $\mathcal{B}^*$ already, i.e. $\beta_{\kappa} = \beta_{\kappa}^*$, and the associated selection probabilities $\{ p_i : i\in \beta_{\kappa}^*\}$ can be correctly determined under $F^*$. In the terminology of Zhang and Oguz-Alper (2020), $\mathcal{B}^*$ is a ``feasible'' BIGS representation of LIS.

%Notice that if it is possible to divide the whole area into grids, where the width of each grid from left to right is determined by the detectability when following a line under LIS, then one can apply second-stage ACS (Thompson, 1991), starting from the grid where a wolverine track intercepts the strip (of grids) selected at the first stage, where each wolverine track runs through a number of contiguous grids that form a network under ACS.

\subsubsection{Estimators}

The HT-estimator is given by Thompson (2012, Ch. 19.1). In the present set-up, since the projection segments are non-overlapping, the selection probability of track $\kappa$ on each draw is given by $p_{(\kappa)} = \sum_{i\in \beta_{\kappa}} p_i$,
where $p_{(1)} = 0.4375$, $p_{(2)} = 0.625$, $p_{(3)} = 0.2$ and $p_{(4)} = 0.5875$. The inclusion probability of $\kappa \in \Omega_s$ is one minus the probability  that $\kappa$ is not selected on any of the 4 draws, $\pi_{(\kappa)} = 1 - (1- p_{(\kappa)})^4$,
where $\pi_{(1)} = 0.90$, $\pi_{(2)} = 0.98$, $\pi_{(3)} = 0.59$ and $\pi_{(4)} = 0.97$. Denote by $p_{(\kappa \cup \ell)} = \sum_{i\in \beta_{\kappa} \cup \beta_{\ell}} p_i$ the probability of selecting either $\kappa$ or $\ell$ on a given draw. The second-order inclusion probability of $\kappa \neq \ell\in \Omega_s$ is given by
\[
%& = 1 - \Big( \mbox{Pr}(\kappa \not\in \Omega_s) + \mbox{Pr}(\ell \not\in \Omega_s) - \mbox{Pr}(\kappa \not\in \Omega_s, \ell \not\in \Omega_s) \Big) \\
\pi_{(\kappa\ell)} = \pi_{(\kappa)} + \pi_{(\ell)} - 1 + \big( 1 - p_{(\kappa \cup \ell)} \big)^4
\]
where $\pi_{(12)} = 0.90$, $\pi_{(13)} = 0.51$, $\pi_{(14)} = 0.88$,  $\pi_{(23)} = 0.57$, $\pi_{(24)} = 0.95$ and $\pi_{(34)} = 0.59$. The HT-estimator $\hat{\theta}_y$ and its estimated variance are given in Table \ref{tab-LIS}.

\begin{table}[ht]
\begin{center}
\caption{IWE under BIGS from $\mathcal{B}^*$.}
\begin{tabular}{l c c c c} \hline
& $\hat{\theta}_y$ & $\hat{\theta}_{z\alpha0}$ & $\hat{\theta}_{z\beta}$ & $\hat{\theta}_{z\alpha.5}$ \\ \hline
Estimate & 7.57 & 9.44 & 8.99 & 9.27 \\
Variance & 5.27 & 1.70 & 2.46 & 1.97 \\ \hline
\end{tabular} \label{tab-LIS}
\end{center}
\end{table}

An unbiased estimator of $\theta$ from the $r$-th draw is $\tau_r = \sum_{\kappa \in \Omega_r} y_{\kappa}/p_{(\kappa)}$, where $\tau_1 = \tau_2 = 7.1878 $ and $\tau_3 = \tau_4 = 11.7021$. Becker (1991) uses the HH-type estimator over all the draws:
\[
\hat{\theta}_{HH} = \sum_{r=1}^4 \tau_r/4 ~.
\]
Under BIGS from $\mathcal{B}^*$, let $\hat{\theta}_{z,r}  = \sum_{i\in s_r} z_i/p_i$ be the IWE on the $r$-th draw, given by \eqref{Zhat}.
Given $\omega_{i\kappa} = p_i/p_{(\kappa)}$ for $i\in \beta_{\kappa}$, where $\omega_{11} = \omega_{43} = \omega_{64} = 1$, $\omega_{12} = p_1/p_{(2)} = 0.7$ and $\omega_{22} = p_2/p_{(2)} = 0.3$, we obtain $\hat{\theta}_{HH}$ above as an IWE estimator, since
\[
\hat{\theta}_{z,r} = \sum_{i\in s_r} \frac{1}{p_i} \sum_{\kappa \in \alpha_i} \frac{p_i}{p_{(\kappa)}} y_{\kappa}
= \sum_{\kappa \in \Omega_r} \frac{y_{\kappa}}{p_{(\kappa)}} = \tau_r ~.
\]
Since these weights are given by \eqref{paid} with $\gamma =0$, we denote $\hat{\theta}_{HH}$ by $\hat{\theta}_{z\alpha0}$ in Table \ref{tab-LIS}.

For the multiplicity estimator $\hat{\theta}_{z\beta}$ with equal weights, we have $\omega_{11} = \omega_{43} = \omega_{64} = 1$, and $\omega_{12} = \omega_{22}= 0.5$. The resulting IWE on each draw are $\hat{\theta}_{z\beta,1} = \hat{\theta}_{z\beta,2} = 6.2736$ and $\hat{\theta}_{z\beta,3} = \hat{\theta}_{z\beta,4} = 11.7021$. The IWE over all the draws is given in Table \ref{tab-LIS}.
Next, the unequal weights by \eqref{paid} can be calculated, since $\alpha_i$ is observed under ancestral BIGS from $\mathcal{B}^*$ given the observation procedure of LIS. Let $\gamma = 0.5$. We have $\omega_{11} = \omega_{43} = \omega_{64} = 1$, $\omega_{12} = 0.6226$ and $\omega_{22} = 0.3773$. The corresponding IWE is 6.8341 on the first two draws, and 11.7021 on the last two draws. This estimator is denoted by $\hat{\theta}_{z\alpha.5}$ in Table \ref{tab-LIS}.

Given the systematic sampling design of the transect lines, the tracks $\{ 1, 2, 4\}$ can only be observed if a position is selected in the left part of 1st projection segment, which would only result in $\{ 1, 5, 6\}$ as the sampled projection segments. Similarly, the tracks $\{ 3,4\}$ can only be observed if a position is selected in 4th projection segment, which would only result in $\{ 4, 6, 7\}$ as the sampled projection segments. Thus, applying the RB method would not change any unbiased IWE based on the observed sample BIGS in this case.

The estimator $\hat{\theta}_{HH}$ of Becker (1991) is the IWE $\hat{\theta}_{z\alpha0}$. The HT-estimator $\hat{\theta}_y$ noted by Thompson (2012) can be given as the IWE with weights satisfying \eqref{Yhat}. Other unbiased IWE can be used for LIS under BIGS from $\mathcal{B}^*$, two of which are as given in Table \ref{tab-LIS}. Neither the HT-estimator $\hat{\theta}_y$ nor the multiplicity estimator $\hat{\theta}_{z\beta}$ is efficient here. Efficiency gains can be achieved using the PIDA weights \eqref{paid}. In this case, adjusting the equal weights by the selection probability while disregarding the degrees of the initial sample units performs well, where $\hat{\theta}_{z\alpha0}$ has the lowest estimated variance. Of course, the true variance of $\hat{\theta}_{z\alpha0}$ may or may not be smaller than that of, say, $\hat{\theta}_{z\alpha.5}$. Meanwhile, setting $\gamma = 1.227$ would numerically reproduce the equal weights $\omega_{12} = \omega_{22}  = 0.5$ based on the observed sample.  It seems that the IWE by \eqref{paid} has the potential to approximate the relatively more efficient estimators in different situations, if one is able to choose the coefficient $\gamma$ in \eqref{paid} appropriately.

\subsection{A simulation study}

\begin{figure}[h]
\hspace{3mm}
\includegraphics[scale=0.35]{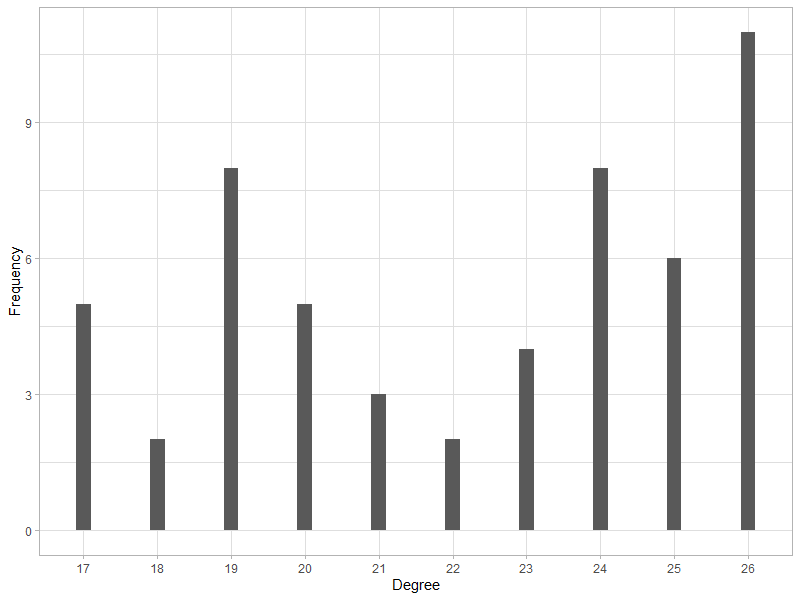}
\includegraphics[scale=0.35]{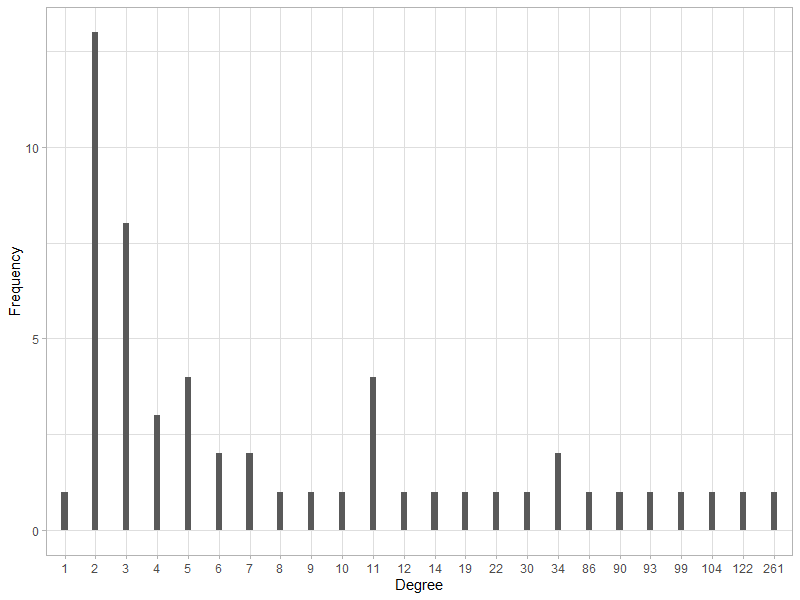}
\caption{Distribution of degree $|\alpha_i|$ in $\mathcal{B}$ (left) and $\mathcal{B}'$ (right).} \label{fig:degree}
\end{figure}

Two graphs $\mathcal{B} = (F, \Omega; H)$ and $\mathcal{B}' = (F, \Omega; H')$ are constructed for this simulation study, which have the same $F$ and $\Omega$, where $|F| = 54$ and $|\Omega| = 310$. The two different edge sets have the same number of edges, where $ |H| = |H'| = 1200$, but different distributions of the degree $|\alpha_i|$ over $F$, as shown in Figure \ref{fig:degree}. The distribution is relatively uniform over a small range of values in $\mathcal{B}$, but much more skewed and asymmetric in $\mathcal{B}'$.

Let $\theta = |\Omega|$, or $y_{\kappa} \equiv 1$ for $\kappa\in \Omega$. We consider the following 7 estimators of $\theta$ under BIGS from $\mathcal{B}$ or $\mathcal{B}'$ with SRS of $s$, where $m = |s|$ varies from 2 to 53:
\begin{itemize}[leftmargin=6mm]
\item the IWE $\hat{\theta}_y$ with weights satisfying \eqref{Yhat} which is the HT estimator;
\item the IWE $\hat{\theta}_{z\alpha\gamma}$ with weights \eqref{paid} for $\gamma = 0, 1, 2$, where $\hat{\theta}_{z\alpha0}$ is the multiplicity estimator;
\item the IWE $\hat{\theta}_p$ by \eqref{Phat} which is the priority-rule estimator of Birnbaum and Sirken (1965), where $F$ is arranged in random, ascending or descending order by $|\alpha_i|$, yielding three estimators, denoted by $\hat{\theta}_{pR}$, $\hat{\theta}_{pA}$ and $\hat{\theta}_{pD}$, respectively.
\end{itemize}

Table \ref{tab-sim} gives the relative efficiency of the 6 other estimators against the HT-estimator, for a selected set of initial sample sizes, each based on 10000 simulations of BIGS from either $\mathcal{B}$ or $\mathcal{B}'$. All the results are significant with respect to the simulation error.

\begin{table}[ht]
\begin{center}
\caption{Relative efficiency of IWE (against $\hat{\theta}_y$) for $\mathcal{B}$ and $\mathcal{B}'$, 10000 simulations.}
\begin{tabular}{c | cccccc | cccccc} \hline
& \multicolumn{6}{|c|}{$\mathcal{B}$} & & \multicolumn{5}{c}{$\mathcal{B}'$} \\
m & $\hat{\theta}_{z\alpha0}$ & $\hat{\theta}_{z\alpha1}$ & $\hat{\theta}_{z\alpha2}$ & $\hat{\theta}_{pR}$ & $\hat{\theta}_{pA}$ & $\hat{\theta}_{pD}$ & $\hat{\theta}_{z\alpha0}$ & $\hat{\theta}_{z\alpha1}$ & $\hat{\theta}_{z\alpha2}$ & $\hat{\theta}_{pR}$ & $\hat{\theta}_{pA}$ & $\hat{\theta}_{pD}$ \\ \hline
5 & 0.96 & 0.55 & 0.49 & 0.80 & 1.43 & 0.68  & 1.22 & 0.23 & 0.18 & 0.97 & 1.16 & 0.83 \\
11 & 0.95 & 0.55 & 0.48 & 0.97 & 2.57 & 0.84 & 1.74 & 0.33 & 0.25 & 0.89 & 1.54 & 0.45 \\
17 & 0.99 & 0.57 & 0.51 & 2.34 & 4.98 & 2.57  & 2.67 & 0.51 & 0.39 & 0.82 & 2.30 & 0.24 \\
29 & 1.31 & 0.75 & 0.67 & 26.7 & 30.1 & 33.2 & 7.96 & 1.54 & 1.17 & 12.0 &  12.1 & 29.3 \\ \hline
%35 & 1.64 & 0.95 & 0.84 & 167.2 & 162.8 & 191.3 & 15.9 & 3.08 & 2.34 & 104.0 & 59.2 & 252.7 \\ 
\end{tabular} \label{tab-sim}
\end{center}
\end{table}

Regarding the priority-rule estimator, we notice that all the three estimators $\hat{\theta}_{pR}$, $\hat{\theta}_{pA}$ and $\hat{\theta}_{pD}$ become biased given large enough initial sample size $m$, which happens at $m=45$ for $\mathcal{B}$ where the maximum degree $|\beta_{\kappa}|$ is 10 over $\Omega$, and $m=46$ for $\mathcal{B}'$ where the maximum degree $|\beta_{\kappa}|$ is 9. Moreover, although the variance of any $\hat{\theta}_p$ initially decreases as $m$ increases, the variance starts to increase with $m$ once the latter is larger than a threshold value, somewhere between 10 and 30 in these simulations, so that the performance of $\hat{\theta}_p$ can deteriorate as the initial sample size increases long before it becomes biased.

The sampling variance of $\hat{\theta}_p$ is also affected by the ordering of the sampling units in $F$. The variance tends to be lowest when $F$ is arranged in descending ordering by $|\alpha_i|$, as long as the variance is decreasing with $m$, whereas ascending ordering tends to yield the largest variance. Without prioritisation, the value $z_i$ is a constant of sampling given $\omega_{i\kappa}$. Due the randomness induced by the priority-rule, $z_i$ varies over different samples. A sampling unit with large $|\alpha_i|$ has a large range of possible $z_i$ values. Placing such a unit towards the end of the ordering tends to increase the sample variance of $\{ z_i : i\in s\}$ due to prioritisation, compared to when the same unit is placed towards the beginning of the ordering, because the rule of Birnbaum and Sirken (1965) prioritises the sampling unit in front of the other ancestors. This is a reason why descending ordering by $|\alpha_i|$ may work better than ascending ordering. However, one may not know $\{ |\alpha_i| : i\in F\}$ in practice, in which case applying $\hat{\theta}_p$ given whichever ordering of $F$ can be a haphazard business.

Given initial SRS, the different HH-type estimators here differ only with respect to the use of $|\alpha_i|$ in the PIDA weights \eqref{paid} via the choice of $\gamma$. The equal-weights esmator $\hat{\theta}_{z\alpha0}$ is the least efficient of the three HH-type estimators, especially for $\mathcal{B}'$ where the distribution of $|\alpha_i|$ is more skewed. The differences between the other two estimators $\hat{\theta}_{z\alpha1}$ and $\hat{\theta}_{z\alpha2}$ are relatively small, compared to their differences to $\hat{\theta}_{z\alpha0}$, so that a non-optimal choice of $\gamma \neq 0$ is less critical than simply setting $\gamma =0$. Taken together, these results suggest that the extra effort that may be required to obtain $|\alpha_i|$ is worth considering in practice, and a sensible choice of $\gamma$ depending on the distribution of $|\alpha_i|$ over $F$ if it is known, or $\mathcal{B}_s$ if it is only observed in the sample BIG, is an interesting question to be studied.

Finally, both $\hat{\theta}_{z\alpha1}$ and $\hat{\theta}_{z\alpha2}$ are more efficient than the HT-estimator when $m$ is small, whereas the HT-estimator improves more quickly as $m$ becomes larger, especially for $\mathcal{B}'$. The matter depends on the sampling fractions $|\Omega_s|/|\Omega|$ and $|s|/|F|$, as well as the respective inclusion probabilities of motifs and sampling units. The interplay between them is complex as it depends on the population BIG. Further research is needed in this respect.

\section{Concluding remarks} \label{remark}

In this paper we develop a large class of incidence weighting estimators \eqref{IWE} under BIGS. The IWE is applicable to all situations of unconventional sampling techniques that require a specific observation procedure in addition to an initial sample, which can be represented by ancestral BIGS, including indirect, network, adaptive cluster and list-intercept sampling. The condition \eqref{unbias} ensures exactly design-unbiased IWE, which synthesises and generalises the conditions underlying the other unbiased estimators known in the literature.

The classic HT-estimator from finite-population sampling is shown to be a special case of IWE, with any sample dependent weights satisfying the restriction \eqref{Yhat}, which provides a novel insight. A more general restriction \eqref{W} is given for sample dependent weights. It will be intriguing to investigate other HT-type estimators satisfying this restriction.

The priority-rule estimator invented by Birnbaum and Sirken (1965) is another a special case of IWE. However, it may become biased as the initial sample size increases and behave erratically long before that, such that its application may be a haphazard business if one is unable to control the interplay between the ordering of sampling units and the priority-rule of Birnbaum and Sirken (1965). It remains to be seen whether one is able to overcome these shortcomings by future developments.

The HH-type estimators used in the literature are also members of the proposed class. While it is in principle possible to apply the Rao-Blackwell method to an HH-type estimator to improve its efficiency, the computation may be intractable if the conditional sample space of $s$ is large and/or if the initial sampling design $p(s)$ is not fully specified. However, consideration of the Rao-Blackwell method and the degrees (in the BIG) of the sampling units points to the PIDA weights \eqref{paid} for IWE, as a general alternative to the commonly used equal weights and the corresponding multiplicity estimator. The numerical illustration of line-intercept sampling and the simulation results suggest that the PIDA weights can easily outperform the equal weights. Further study is warranted, in order to identify the sensible choice of the PIDA weights in applications.

Finally, other incidence weights can be explored subjected to the condition \eqref{unbias}, beyond those examined in this paper. This is clearly another direction of future research.

\appendix
\section{Priority probabilities of $\hat{\theta}_p$} \label{pprob}

For each $\kappa \in \Omega_s$ and $i\in\beta_{\kappa}$, let $d_{i(\kappa)} = \sum_{j \in F: j < i} I_{(j\kappa) \in H}$ be the number of sampling units where higher priority than $i$ under the priority rule $\min(s\cap \beta_{\kappa})$. Assume SRS of $s$, where $m= |s|$. We have
\[
p_{i\kappa} = \binom{M-1- d_{i(\kappa)}}{m-1} / \binom{M-1}{m-1} ~.
\]
The joint priority probability of $(i\kappa)$ and $(j\ell)$ given $\delta_i\delta_j =1$ is
\[
p_{i\kappa, j\ell} = \begin{cases} p_{i\kappa} & \text{if } \kappa = \ell, i = j \\
0 & \text{if } \kappa = \ell, i \neq j \\
\binom{M-1- d_{i(\kappa, \ell)}}{m-1} / \binom{M-1}{m-1} & \text{if } \kappa \neq \ell, i = j \\
\binom{M-2- d_{i(\kappa), j(\ell)}}{m-2} / \binom{M-2}{m-2} & \text{if } \kappa \neq \ell, i \neq j \text{ and }
|\beta^i_{\kappa} \cap \{ j \} | + |\beta^j_{\ell} \cap \{ i\} | = 0 \\
0 & \text{if } \kappa \neq \ell, i \neq j \text{ and }  |\beta^i_{\kappa} \cap \{ j\} | + |\beta^j_{\ell} \cap \{ i\} | > 0
\end{cases}
\]
where $\beta^i_{\kappa}$ is the subset ancestors of $\kappa$ with higher priority than $i$, and $d_{i(\kappa, \ell)} = | \beta^i_{\kappa} \cup \beta^i_{\ell} |$ is the number of units in $\beta_{\kappa} \cup \beta_{\ell}$ with higher priority than $i$, and $d_{i(\kappa), j(\ell)} = | \beta^i_{\kappa} \cup \beta^j_{\ell} |$.

\end{document}